\documentclass[11pt]{amsart}

\usepackage{color}
\usepackage{amsfonts}
\usepackage{amscd}
\usepackage{amsmath}
\usepackage{amssymb}
\usepackage{dsfont}
\usepackage{mathabx}

\usepackage{tikz}
\usetikzlibrary{matrix,arrows,decorations.pathmorphing}
\usepackage{pgfkeys}
\newtheorem{Th}{Theorem}[section]
\newtheorem{lemma}[Th]{Lemma}
\newtheorem{prop}[Th]{Proposition}
\newtheorem{defin}[Th]{Definition}
\newtheorem{cor}[Th]{Corollary}

\newcommand{\Proof}{\noindent {\it Proof. }}
\newcommand{\dist}{\mathop{\rm dist}}

\newcommand{\Ker}{\mathop{\rm Ker}}
\newcommand{\Ran}{\mathop{\rm Ran}}

\newcommand{\lra}{\longrightarrow}
\newcommand{\Lra}{\Leftrightarrow}

\newcommand{\A}{\mathcal{A}}

\newcommand{\Lc}{\mathcal{B}}
\newcommand{\Cl}{\mathcal{C}}
\newcommand{\D}{\mathcal{D}}

\newcommand{\Ss}{\mathcal{SS}}
\newcommand{\SC}{\mathcal{SC}}

\newcommand{\K}{\mathcal{K}}

\newcommand{\W}{\mathcal{W}}

\newcommand{\N}{\mathbb{N}}

\newcommand{\Rr}{\mathbf{R}}
\newcommand{\Sg}{\mathcal{S}}
\newcommand{\Cg}{\mathcal{C}}

\newcommand{\tendsstyle}[5]{\relax\ifmmode{%
 \setbox0\hbox{$#1#3$}\dimen0 \ht0 \advance\dimen0 -\fontdimen22#22
 #1\mathop{\vcenter to 2\dimen0{\box0\vss}}\limits_{#4}^{#5}%
}\fi}

\newcommand{\fin}{\hspace*{\fill} $\square$\vspace{0.5\baselineskip}}

\begin{document}

\title[Tauberian pairs of subspaces]
  {Tauberian pairs of closed subspaces \\ of a Banach space}

\thanks{Supported in part by Gobierno de Cantabria. Contrato Programa Gob. Cantabria-UC 2024\newline
}

\author[M. Gonz\'alez]{Manuel Gonz\'alez}
\address{Departamento de Matem\'aticas, Facultad de Ciencias, Universidad de Cantabria, E-39071 Santander, Espa\~na}
\email{manuel.gonzalez@unican.es}

\author[A. Mart\'\i nez-Abej\'on]{Antonio Mart\'\i nez-Abej\'on}
\address{Departamento de Matem\'aticas, Facultad de Ciencias, Universidad de Oviedo, E-33007 Oviedo, Espa\~na}
\email{ama@uniovi.es}
%
%\author[J. Pello]{Javier Pello}
%\address{Escuela Superior de Ciencias Experimentales y Tecnolog\'\i a, Universidad Rey Juan Carlos, 28933 M\'ostoles, Spain}
%\email{javier.pello@urjc.es}

%\date{}
%\date{\today}

\dedicatory{To the Memory of Albrecht Pietsch 1934-2024}
%%% ----------------------------------------------------------------------

\begin{abstract}
We introduce the notions of tauberian, cotauberian and weakly compact pair of closed subspaces of  a Banach space.
The theory produced by these notions is richer than that of the corresponding operators since an operator can be regarded as a suitable pair of closed subspaces.
We investigate into these classes of pairs of subspaces and describe several applications in order to define some notions of indecomposability for Banach spaces and
 in order to extend definitions from the case of bounded operators to the case of closed operators.
\end{abstract}

\maketitle
\thispagestyle{empty}

%\tableofcontents

\section{Introduction}

Arguably, Pietsch's notion of operator ideal \cite{Pietsch:80} was inspired by the class of compact operators in Fredholm theory. Similarly, the semi-Fredholm operators led to the notion of \emph{operator semigroup} in \cite{AGM:01}: for every operator ideal $\A$ we have two operator semigroups $\A_+$ and $\A_-$ defined by
\begin{equation}\label{A+-}
\A_+:=\{T\colon TS\in\A\Rightarrow S\in \A\} \textrm{ and } \A_-:=\{T\colon ST\in\A\Rightarrow S\in \A\}.
\end{equation}

For the ideal $\K$ of compact operators, $\K_+$ and $\K_-$ are, respectively, the classes of upper semi-Fredholm and lower semi-Fredholm operators, and for the ideal $\W$ of weakly compact operators, $\W_+$  is the class of tauberian operators, introduced by Kalton and Wilansky \cite{KW:76}, and  $\W_-$ is the class of cotauberian operators introduced by Tacon \cite{Tac:83}.

In \cite[Section IV.4]{Kato:66}, Kato introduced a Fredholm theory for pairs of closed subspaces of a Banach space that allowed him to extend some small-norm perturbation results for bounded operators to the case of closed operators, and the notions of compact, strictly singular and strictly cosingular pairs of subspaces were introduced in \cite{Gonzalez:05}.

Here we extend the concepts of tauberian, cotauberian and weakly compact operators to pairs of closed subspaces of a Banach space. These  operators provided a  meaningful Fredholm theory in which compact operators and finite dimensional spaces are replaced by weakly compact operators and reflexive spaces, and they have found many applications in Banach space theory (see \cite{GM:10}).
Observe that there are many examples of tauberian and cotauberian operators since every operator $T$ can be factorized as $T= B A$ where $A$ is cotauberian and $B$ is tauberian (see \cite{DFJP} and \cite[Section 3.2]{GM:10}).
\medskip

Let us describe the contents of the paper. In Section \ref{sect:prelim}, we give some definitions and results that connect bounded and closed operators acting between Banach spaces with pairs of closed subspaces of a Banach space. Moreover, given a class  of bounded operators invariant under bijective isomorphisms, we define the pairs of closed subspaces in $\mathcal{A}$, and we introduce some properties for classes of pairs.

In Section \ref{sect:tau-pairs} we recall some characterizations of tauberian, cotauberian and weakly compact operators, and apply them to characterize pairs of closed subspaces $(M,N)$ of a Banach space $Z$ in the corresponding classes and give some properties of these classes of pairs. We also show that any weakly compact pair $(M,N)$ remains weakly compact if one of its entries is perturbed by a small weakly compact perturbation (Proposition \ref{useful-prop}), that  tauberian pairs are preserved in certain sense by tauberian operators and cotauberian pairs are preserved by cotauberian operators (Theorem \ref{th:preserve}), and we give two versions for pairs (Theorems \ref{pert-char} and \ref{pert-char-alt}) of the perturbative characterization of tauberian and cotauberian operators (parts (4) of Propositions \ref{tau-op} and \ref{cotau-op}).

In Section \ref{sect:appl} we apply the previous results to describe some notions of indecomposability for Banach spaces defined in terms of tauberian and cotauberian pairs, we show how to extend the notions of weakly compact, tauberian and cotauberian bounded operators to closed operators, giving an alternative version to the concepts  studied by Cross in \cite{Cross:92,Cross:95} and \cite{Cross:98}, and we extend the notions of  $\Rr\Sg\Sg$ and $\Rr\Sg\Cg$ operators of Hermann \cite{Herman:68} to pairs of closed subspaces. We also pose some open questions.
\medskip

Notation: $X$, $Y$ and $Z$ are Banach spaces, $B_X$ and $S_X$ are the closed unit ball and the unit sphere of $X$, and $\Lc(X,Y)$ is the set of all bounded operators from $X$ into $Y$.  We say that $T\in\Lc(X,Y)$ is an isomorphism if it has a continuous inverse on its range.

If $\mathcal{O}$ is a class of bounded operators, its components are $\mathcal{O}(X,Y)=\mathcal{O}\cap\Lc(X,Y)$.
The class $\mathcal{O}$ is \emph{stable by bijective isomorphisms} if $R S T\in \mathcal{O}\cap\Lc(W,Z)$ for all $S\in\mathcal{O}(X,Y)$ and all bijective isomorphisms $T\in\Lc(W,X)$ and $R\in\Lc(Y,Z)$; the class $\mathcal{O}$ is \emph{open} if for every $X$ and $Y$, $\mathcal{O}(X,Y)$ is open in $\Lc(X,Y)$; given a second class $\mathcal{D}$
of bounded operators, $\mathcal{O}+\mathcal{D}\subset\mathcal{O}$ means that given $S\in\mathcal{O}(X,Y)$ and $K\in\mathcal{D}(X,Y)$, $S+K\in\mathcal{O}$.

For an operator $T:D(T)\subset X\to Y$ with \emph{domain} $D(T)$, the \emph{graph} is $G(T)=\{(x,Tx) : x\in D(T)\}$. We say that $T$ is \emph{closed} if $G(T)$ is a closed subspace of $X\times Y$. Also, $\Cl(X,Y)$ denotes the set of all closed operators with domain in $X$ and range contained in $Y$. Note that $\Lc(X,Y)$ is a subset
of $\Cl(X,Y)$. For short, $\Lc(X)$ and  $\Cl(X)$ respectively denote $\Lc(X,X)$ and $\Cl(X,X)$. For any $T\in \Cl(X,X)$, we respectively denote $\Ker T$ and $\Ran T$ the kernel of $T$ and the range of $T$.

In addition, $X^*$ is the dual space of $X$, and the annihilator of a subspace $A$ of $X$ is $A^\perp= \{f\in X^* \colon f(a)=0\; \textrm{for all $ a\in A$}\}$. Given $S\in\Lc(X,Y)$, $S^*\in\Lc(Y^*,X^*)$ is the conjugate operator of $S$. As usual, $X$ is identified with a subspace of its second dual $X^{**}$.  Thus, the second conjugate $S^{**} \in\Lc(X^{**} ,Y^{**})$ is an extension of $S$.
%For $T\in C_D(X,Y)$, the conjugate operator is defined as follows: $D(T^*)= \{g\in Y^* \colon \hbox{$g\circ T$ is continuous}\}$, and for $g\in D(T^*)$, $T^*g\in X^*$ is the (unique) continuous extension of $g\circ T$.
%
We denote by $\mathcal{F}_Z$ the set of all ordered pairs of closed subspaces of $Z$. For any closed subspace $M$ of $Z$,  $J_M$ is the inclusion operator of $M$ into $Z$ and $Q_M$ is the quotient map from $Z$ onto $Z/M$.

\section{Preliminaries}\label{sect:prelim}

Sometimes in this paper, if $A$ is a subspace of $X$ and $B$ is a subspace of $Y$, we identify $A$ with $A\times \{0\}$
and $B$ with $\{0\} \times B$, both subspaces of $X\times Y$.

Every $(M,N)\in\mathcal{F}_Z$ determines an operator $Q_NJ_M\in\Lc(M,Z/N)$.
Conversely, every $S\in\Lc(X,Y)$ is determined up to bijective isomorphism by each of the pairs of closed subspaces
$(G(S),X)$ and $(X,G(S))$ of $X\times Y$. Indeed,
\begin{equation}\label{rep-op}
 S=Q_XJ_{G(S)}U \quad \textrm{ and } S=VQ_{G(S)}J_X,
\end{equation}
where $U\colon X\lra G(S)$ and $V\colon (X\times Y)/G(S)\lra Y$ are bijective isomorphisms given by
$Ux=(x,Sx)$ and $V((x,y)+G(S))= Sx-y$ \cite[Proposition 2.5]{Gonzalez:05}, and $Q_X\colon X\times Y\lra (X\times Y)/(X\times\{0\})\equiv Y$ is the quotient operator that maps $(x,y)$ to $y$.
\smallskip

Kato calls \emph{Fredholm} to any pair $(M,N)\in\mathcal{F}_Z$ for which $Q_NJ_M$ is a Fredholm operator (Chapter IV,\S 4 in \cite{Kato:66}). Following the same spirit,
we set the following definition which associates each class of operators stable by bijective isomorphisms with a class of pairs of subspaces:

\begin{defin}\label{def:pair-in-O}
If $\mathcal{O}$ is a class of bounded operators stable by bijective isomorphisms, we will denote by the same symbol $\mathcal{O}$ the collection of all pairs $(M,N)\in\mathcal{F}_Z$ for any Banach space $Z$ for which $Q_NJ_M\in\mathcal{O}$.
\end{defin}

The reader should beware of the difference between $T\in\mathcal{O}$ and $(M,N)\in\mathcal{O}$.

The next statements hold true for any $T\in \Cl(X,Y)$. Their proof is a simple exercise.

\begin{enumerate}
\item[(i)] $G(T)\cap X = \{x\in D(T): Tx=0\} = \Ker T$ is a closed subspace of $X$.
\smallskip

\item[(ii)] $X+G(T) = \{(x+z,Tz) : x\in X, z\in D(T)\} = X\times \Ran T$.

Thus $\Ran T$ is closed if and only if $X+G(T)$ is closed, and
$$
\dim Y/\Ran T=\dim (X\times Y)/(X+G(T)).
$$
%\item[(iii)] If $D(T)$ is dense then $G(T^*)=\{(Tx,-x) :x\in D(T)\}^\perp$ \cite{Goldberg:66}. In particular, the conjugate operator is closed: $T^*\in \Cl(Y^*,X^*)$.
\end{enumerate}
%\smallskip

The definitions of bounded semi-Fredholm operators are easily extended to the context of closed operators as follows: $T\in\Cl(X,Y)$ is said to be \emph{upper semi-Fredholm} (resp. \emph{lower semi-Fredholm}) if $\Ker{T}$ is finite dimensional and $\Ran{T}$ is closed ($\Ran{T}$ is closed and finite codimensional in $Y$).
A translation into the language of our Definition \ref{def:pair-in-O} asserts that  $T\in\Cl(X,Y)$ is upper (lower) semi-Fredholm if and only if $(G(T),X)\in\K_+$ (resp. $(G(T),X)\in\K_-$), %\cite[Chapter IV]{Kato:66}
a fact that can be directly verified by the reader from (i) and (ii). Therefore, Definition \ref{def:pair-in-O} also builds a link  between classes of pairs of subspaces and classes of closed operators.
%This fact will be used in Section \ref{sect:appl}.
\medskip

Next we state some properties of the kernel and  the range  of the operator $Q_NJ_M$ associated with the pair $(M,N)\in\mathcal{F}_Z$:
\begin{enumerate}
\item[(i)] $\Ran(Q_NJ_M)=(M+N)/N$ and $\Ker(Q_NJ_M)=N\cap M$.
\item[(ii)] $\Ran(Q_NJ_M)$ closed $\Lra M+N$ closed $\Lra M^\perp+N^\perp$ closed \cite[Chapter IV, \S 4,2]{Kato:66}.
\item[(iii)] $(Q_NJ_M)^* = Q_{M^\perp}J_{N^\perp}$.
\item[(iv)] The dual of $(X/N)/\overline{\Ran Q_NJ_M}$ is (linearly isometric to) $\Ker Q_{M^\perp} J_{N^\perp}$.
\item[(v)] The dual of $\overline{\Ran Q_NJ_M}$ is (linearly isometric to) $Z^*/\Ker Q_{M^\perp} J_{N^\perp}$.
\item[(vi)] If $\Ran Q_NJ_M$ is closed, then the dual of $\Ker Q_NJ_M$ is (linearly isometric to) $Z^*/\Ran Q_{M^\perp} J_{N^\perp}$.
%\item
\end{enumerate}

Let us list some interesting properties that a pair of subspaces may possess:

\begin{defin}\label{defin-alg}
Let $\D$ and $\mathcal{O}$ be two classes of pairs of subspaces.
\begin{itemize}
\item $\D$ is \emph{symmetric} if a pair $(M,N)\in \D$ if and only if $(N,M)\in \D$.
\item $\D$ is \emph{transitive} if both $(L,M), (M,N)\in \D$ implies $(L,N)\in\D$.
\item $\D$ is \emph{left stable under perturbation by $\mathcal{O}$} if $(L,M)\in\mathcal{O}$ and $(M,N)\in\D$ imply $(L,N)\in \D$.
\item $\D$ is \emph{right stable under perturbation by $\mathcal{O}$} if $(L,M)\in\D$ and $(M,N)\in\mathcal{O}$ imply $(L,N)\in \D$.
%\item $\D$ is \emph{stable under perturbation by $\E$} if it is left stable and right stable.
\end{itemize}
\end{defin}

Clearly $\K_+$ and $\K_-$ are symmetric, and in \cite[Theorems 3.1 and 3.9]{Gonzalez:05} it was shown that $\K_+$ is left stable under $\Ss$, $\K_-$ is right stable under $\SC$, and $\K$, $\Ss$, and $\SC$ are transitive. However, $\K$ and $\W$ are not symmetric. Take, for example, a non-reflexive Banach space $Z$, $M=\{0\}$ and $N=Z$. This fact reveals that the classes of pairs of subspaces are richer than the corresponding classes of operators. As we can see in Equation (\ref{rep-op}), an operator $S\in\Lc(X,Y)$ is represented by the pair $(G(S),X)$ and also by its symmetric pair $(X,G(S))$. Hence, these two pairs share the same properties.

\section{Tauberian, cotauberian and weakly compact pairs}\label{sect:tau-pairs}

The original definitions of tauberian and cotauberian operators are the following:

\begin{defin}  (\cite{KW:76}, \cite{Tac:83}) An operator $S\in\Lc(X,Y)$ is said to be:
\begin{enumerate}
  \item[(i)] tauberian if $S^{**}(X^{**}\setminus X)\subset Y^{**}\setminus Y$;
  \item[(ii)] cotauberian if $S^*$ is tauberian.
\end{enumerate}
\end{defin}

It was proved in \cite{GO:90} that the classes of tauberian and cotauberian operators are the operator semigroups $\W_+$ and $\W_-$ associated with the operator ideal $\W$ of weakly compact operators. Both tauberian and cotauberian operators
are studied in \cite{AsTy:90}, \cite{Holub:93}, \cite{JohnsonNST:15}, \cite{NR:85} and \cite{Yang:76}.

We denote by $X^{co}$ the quotient $X^{**}/X$.
The \emph{residuum} of $S\in\Lc(X,Y)$ is the operator $S^{co}:X^{co}\to Y^{co}$ defined by  $S^{co}(x^{**}+X):= S^{**}x^{**} +Y$
(see \cite{GSaTy:95} and \cite[Section 3.1]{GM:10} for the properties of $S^{co}$).

%Next we collect the main  characterizations of the tauberian and the cotauberian operators in order to be used later.

\begin{prop}\label{tau-op} \cite{GO:90,KW:76}
For $S\in\Lc(X,Y)$, the following assertions are equivalent:
\begin{enumerate}
\item[(1)] $S$ is tauberian.
\item[(2)] $S^{co}$ is injective.
\item[(3)] For every compact operator $K\in \K(X,Y)$, $\Ker (S+K)$ is reflexive.
\item[(4)] Every bounded sequence $(x_n)$ in $X$ with $(Sx_n)$ convergent contains a weakly convergent subsequence.
%\item[(6)] $S(B_X)$ closed and $\Ker S=\Ker S^{**}$.
\item[(5)] $\overline{S(B_X)}\subset \Ran S$ and $\Ker S= \Ker S^{**}$.
\end{enumerate}
\end{prop}

\begin{prop}\label{cotau-op} \cite{GO:90,KW:76}
For $S\in\Lc(X,Y)$, the following assertions are equivalent:
\begin{enumerate}
\item[(1)] $S$ is cotauberian.
\item[(2)] $S^{co}$ has dense range.
\item[(3)] For every compact operator $K\in \K(X,Y)$, $Y/\overline{\Ran S+K}$ is reflexive.
\item[(4)] $\Ker S^*=\Ker S^{***}$.
\end{enumerate}
\end{prop}

%\begin{rem}\label{remark}
%Note that an operator $S\in\Lc(X,Y)$ is weakly compact if and only if for every bounded sequence $(x_n)$ in $X$, $(Sx_n)$ has a weakly convergent subsequence; equivalently, %its conjugate operator $S^*\in \W$ (Gantmacher's theorem); or equivalently, $S^{co}=0$.
%\end{rem}

%\begin{cor}\label{ker-coker}
%Let $S\in\Lc(X,Y)$ and $K\in\W(X,Y)$.

%If $S$ is tauberian then so is $S+K$ and $\Ker S$ is reflexive.

%If $S$ is cotauberian then so is $S+K$ and $Y/\overline{\Ran S}$ is reflexive.
%%\end{cor}
%\begin{proof}
%Note that $S^{co}=0 \Rightarrow\Ker S= \Ker S^{**}\Rightarrow \Ker S$ reflexive, and that $(Y/\overline{\Ran S})^*$ can be identified with $\ker S^*$.
%\end{proof}

%%\begin{rem}\label{rem:stau}\;
%\begin{enumerate}
%%\item In general, the sets $\W_+(X,Y)$ and $\W_-(X,Y)$ are not open in $\Lc(X,Y)$ \cite{AlG:91}.
%\item $S^{**}$ tauberian implies $S$ is tauberian. However, the converse implication fails in general \cite{AlG:91}.

%Recall that $S\in\Lc(X,Y)$ is \emph{strongly tauberian} if $S^{co}$ an (into) isomorphism, and $S$ is \emph{strongly cotauberian} if $S^*$ is strongly tauberian %\cite{Rosenthal:99}.

%\item $S$ strongly tauberian if and only if so is $S^{**}$ \cite{Rosenthal:99}, and the sets of strongly tauberian and strongly cotauberian operators are open in the space of %all bounded operators.

%\item In some cases; e.g., when $X=L_1(\mu)$, $S\in\Lc(X,Y)$ is strongly tauberian whenever it is tauberian.
%%\end{enumerate}
%\end{rem}

%\section{Tauberian and cotauberian pairs of closed subspaces}\label{sect:tau-pairs}

For any closed subspace $M$ of $Z$, the pair $(M,M)$ is associated with the null operator $Q_MJ_M=0$. Hence, $(M,M)\in\W$ in a trivial way.
A small perturbation carried on $M$  produces a subspace $N$ so that
 both $(M,N)$ and $(N,M)$ remain  in $\W$ as proved in the following:

\begin{prop}\label{useful-prop}
Let $M$ be a closed subspace of a Banach space $Z$  and let $K\in\W(Z)$ with $\|K\|<1$. Then, for $N=(I-K)(M)$, both $(M,N)$ and $(N,M)$ belong to $\W$.
\end{prop}
\Proof
The operator $U=I-K$ is a bijective isomorphism such that $U(M)=N$. Thus, it induces two bijective  isomorphisms $U_s:M\to N$ and $U_q:Z/M\to Z/N$ that satisfy $UJ_M= J_NU_s$ and $Q_NU=U_qQ_M$. Moreover $(I-K)^{-1} = I+L$ with $L= \sum_{j=1}^\infty K^j\in\W$. We have $(M,N)\in\W$ because
$$
Q_NJ_M= Q_NU^{-1}UJ_M= Q_N(I+L)J_NU_s= Q_NLJ_N U_s\in\W.
$$

Similarly, since $U^{-1}J_N=J_MU_s^{-1}$,
\[
Q_MJ_N= Q_MUU^{-1}J_N= Q_M(I-K)J_MU_s^{-1}= -Q_MKJ_M U_s^{-1} \in\W. \qedhere
\]
\fin

For a subspace $M$ of $Z$, we denote by $M^{co}$ the range of $(J_M)^{co}$, so that  $\Ker ({Q_M}^{co})= M^{co}$ (see  \cite[Proposition 3.1.13]{GM:10}). This enables us to adopt the following identifications:
\begin{equation}
 (J_M)^{co}\equiv J_{M^{co}}\ \hbox{ and } \ (Q_M)^{co}\equiv Q_{M^{co}}. \label{eq:JNco}
\end{equation}
It is also proved in \cite[Corollary 3.1.14]{GM:10} that
\begin{equation}
M^{co}\equiv \frac{M^{\perp\perp}}{M} \equiv \frac{M^{\perp\perp}+Z}{Z}. \label{eq:JNcos}
\end{equation}

Observe that $(M,N)\in \W_+ \Rightarrow M\cap N = M^{\perp\perp}\cap N^{\perp\perp} \Rightarrow M\cap N$ reflexive.

The following three results characterize tauberian pairs, cotauberian pairs and weakly compact pairs. They easily follow from Propositions \ref{tau-op} and \ref{cotau-op} and the basic properties of weakly compact operators.

\begin{prop}\label{tau-pair}
For a pair $(M,N)\in\mathcal{F}_Z$, the following assertions are equivalent:
\begin{enumerate}
\item[(1)] $(M,N)\in \W_+$;
\item[(2)] $M^{co}\cap N^{co}=\{0\}$;
\item[(3)] $M^{\perp\perp}\cap N^{\perp\perp}= M\cap N$ and $\overline{Q_NB_M}\subset (M+N)/N$.
\end{enumerate}
\end{prop}
\Proof
By means of the identifications given in Equation (\ref{eq:JNco}), it follows that $(M,N)\in\W_+$ if and only if
$(Q_NJ_M)^{co}=Q_{N^{co}} J_{M^{co}}$ injective. The remainder is a consequence of Proposition \ref{tau-op}.
 \fin

\begin{prop}\label{cotau-pair}
For a pair $(M,N)\in\mathcal{F}_Z$, the following assertions are equivalent:
\begin{enumerate}
\item[(1)] $(M,N)\in \W_-$;
\item[(2)] $(N^\perp,M^\perp)\in \W_+$;
\item[(3)] $\Ran Q_{N^{co}} J_{M^{co}}$ dense; i.e. $M^{co}+ N^{co}$ is dense in $Z^{co}$.
%\item
%\item
\end{enumerate}
\end{prop}
\Proof
(1)$\Leftrightarrow$(2) By definition, $(M,N)\in\W_-$ if and only if $(Q_NJ_M)^*=Q_{M^{\bot}}J_{N^{\bot}}\in\W_+$, in other words, if and only if $(N^\perp,M^\perp)\in \W_+$.

(2)$\Leftrightarrow$(3) It is a direct consequence of Proposition \ref{cotau-op} and Equation (\ref{eq:JNcos}).
\fin

\begin{prop}  \label{weakcompact-pair}
Let $(M,N)\in\mathcal{F}_Z$. Then  $(M,N)\in \W$ if and only if  $(N^\perp,M^\perp)\in \W$; or equivalently,  $M^{co}\subset N^{co}$.
\end{prop}
\Proof
The proof follows from the fact that $S\in\W$ if and only if $S^*\in\W$ (Gantmacher theorem) and the fact that $S\in\W$ if and only if $S^{co}=0$.
\fin

The proof of the next Lemma is easy.

\begin{lemma}\label{lemma}
If $T\in\Lc(X,Y)$ and $M$ is a closed subspace of $X$ then $\overline{T(M)}^\perp = T^{*-1}(M^\perp)$.
\end{lemma}
%\begin{proof}
%For the direct inclusion, we take $g\in T(M)^\perp$ and show that $T^*g\in M^\perp$. Indeed, for every $m\in M$, $\langle m,T^*g \rangle= \langle Tm,g\rangle= 0$.

%For the converse, we take $g\in Y^*$ such that $T^*g \in M^\perp$ and show that $g\in \overline{T(M)}^\perp$. Indeed, if $y=Tm$ with $m\in M$ then $\langle y,g \rangle= \langle Tm,g\rangle= \langle m,T^*g\rangle=0$.
%\end{proof}

Propositions \ref{tau-pair} and \ref{cotau-pair} lead us to some stronger results about the preservation of tauberian (cotauberian) pairs by tauberian (cotauberian) operators:

\begin{Th}\label{th:preserve}
Let $T\in \Lc(X,Y)$.
\begin{enumerate}
\item[(i)] If $T$ is tauberian and $(M,N)\in\mathcal{F}_Y$ is tauberian then $(T^{-1}(M),T^{-1}(N))$ is tauberian.
\item[(ii)] If $T$ is cotauberian and $(M,N)\in\mathcal{F}_X$ is cotauberian then $(\overline{T(M)}, \overline{T(N)})$ is cotauberian.
\end{enumerate}
\end{Th}
\Proof
(i) Let us denote $A=T^{-1}(M)$ and $B=T^{-1}(N)$. In view of Proposition \ref{tau-pair}, we have to show that $A^{co}\cap B^{co}=0$. Since $T^{co}$ is injective and $M^{co}\cap N^{co}=\{0\}$, we have $T^{co-1}(M^{co})\cap T^{co-1}(N^{co})=\{0\}$. Thus, it is enough to show that $A^{co}\subset T^{co-1}(M^{co})$; equivalently, $w\in A^{co}$ implies $T^{co}w\in M^{co}$. Indeed, if $w\in A^{co}$ then $w=x^{**}+X$ with $x^{**}\in A^{\perp\perp}$. Moreover,  $T(A)\subset M$ implies $T^{**}(A^{\perp\perp})\subset M^{\perp\perp}$. Thus $T^{co}w= T^{**} x^{**}+Y\in (M^{\perp\perp}+ Y)/Y\equiv M^{co}$, and the proof of this part is done.

(ii) In this case $T^*\in\Lc(Y^*,X^*)$ is tauberian and $(N^\perp,M^\perp)$ is tauberian in $X^*$. Thus, part (i) shows that $(T^{*-1}(N^\perp), T^{*-1}(N^\perp))$ is tauberian in $Y^*$ and, by Proposition \ref{cotau-pair} and Lemma \ref{lemma}, $(\overline{T(M)}, \overline{T(N)})$ is cotauberian.
\fin

The following corollary is a straightforward consequence from the fact that all isomorphisms are tauberian operators and all surjections are cotauberian.

\begin{cor}
(1) If $Z_1$ is a closed subspace of $Z_2$ and $M,N\in\mathcal{C}(Z_1)$  then $(M,N)$ tauberian in $Z_2$ implies $(M,N)$ tauberian in $Z_1$.

(2) If $Z_2$ is a quotient of $Z_1$, $Q:Z_1\to Z_2$ is the quotient map and $M,N\in \mathcal{C}(Z_1)$ then $(M,N)$ cotauberian in $Z_1$ implies $(\overline{Q(M)}, \overline{Q(N)})$ cotauberian in $Z_2$.
\end{cor}

Next, we show some properties of the classes of pairs $\W$, $\W_+$ and $\W_-$:

\begin{prop}\label{prop:basic}\;
\begin{enumerate}
\item[(i)] The classes $\W_+$ and $\W_-$ are symmetric.
\item[(ii)] The class $\W$ is transitive but non-symmetric.
\item[(iii)] $\W_+$ is left stable but not right stable under perturbation by $\W$.
\item[(iv)] $\W_-$ is right stable but not left stable under perturbation by $\W$.
\end{enumerate}
\end{prop}
\Proof
(i) By Proposition \ref{tau-pair}, the pair $(M,N)$ is tauberian if and only if $M^{co}\cap N^{co}= \{0\}$, which is a symmetric property. The result for $\W_-$ follows by duality.

\emph{An alternative proof for the symmetry of $\W_+$:}
Recall that $T\in\W_+$ if and only if $\Ker T^{**}= \Ker T$ and $\overline{T(B_X)}\subset \Ran T$ (part (4) in Proposition \ref{tau-op}).
Since $M\cap N = M^{\perp\perp}\cap N^{\perp\perp}$ is symmetric, it is enough to show that $\overline{Q_N(B_M)} \subset (M+N)/N$ implies $\overline{Q_M(B_N)}\subset (M+N)/M$.

To do that, let $(v_n)\subset B_N$ such that $(Q_M v_n)$ converges to $Q_M x$. We have to show that $x$  can be chosen in $N$. Indeed, since $\dist(v_n-x,M)\to 0$ there exists $(u_n)\subset M$ such that $\|u_n+v_n-x\|\to 0$, and $Q_N(u_n+v_n-x) = Q_N(u_n-x)\to 0$. Thus, $(u_n)$ bounded and $Q_N(B_M)\subset (M+N)/N$ imply $Q_N x \in (M+N)/M$, and we can assume $x\in N$.
\medskip

(ii) Note that $(L,M)\in\W\Lra (Q_MJ_L)^{co}= Q_{M^{co}} J_{L^{co}}=0\Lra L^{co}\subset M^{co}$. So it is enough to observe that $L^{co}\subset M^{co}, M^{co}\subset N^{co} \Rightarrow L^{co}\subset N^{co}$.

To see that $\W$ is not symmetric, let $Z$ be a non-reflexive Banach space. Then $(Z,\{0\})$ corresponds to the identity on $Z$, and $(\{0\},Z)$ corresponds to the operator $0$.

(iii) Arguing as is (1) and (2), $L^{co} \subset M^{co}$ and $M^{co}\cap N^{co}= \{0\}$ imply $L^{co} \cap N^{co}= \{0\}$, proving the first part.

For the second part, take $L$ non-reflexive, $M=\{0\}$ and $N=L$. Then $(L,\{0\})\in \W_+$ and $(\{0\},L) \in\W$, but $(L,L)\notin \W_+$.

(iv) If $(L,M)\in \W_-$ and $(M,N)\in\W$ then
$(N^\perp,M^\perp)\in \W$ and $(M^\perp, L^\perp) \in\W_+$. By (3), $(N^\perp,L^\perp) \in \W_+$, hence $(L,N)\in\W_-$.

For the second part, take $L$ such that $L^\perp$ is non-reflexive, $M=Z$ and $N=L$. Then $(L,Z)\in \W$ and $(Z,\{0\}) \in\W_-$, but $(L,\{0\})\notin \W_-$.
\fin

The following result corresponds to the perturbative characterizations for tauberian and cotauberian operators, included as part (4) in Propositions \ref{tau-op} and \ref{cotau-op}.
%For example, a pair $(L,N)$ can be regarded as the pair $(M,N)$ after perturbation on the left by the compact pair $(L,M)\in\K$ (see Definition \ref{defin-alg}).
% if $(L,M)\in \K$, we consider $(L,N)$ as a compact perturbation to the left of $(M,N)$, and $L\cap N$ is the kernel of $(L,N)$.

\begin{Th}[Perturbative characterizations]\label{pert-char}
Let $(M,N)\in \mathcal{F}_Z$.
\begin{enumerate}
\item[(i)] The pair $(M,N)$ is tauberian if and only if for every closed subspace of $Z$ with  $(L,M)\in\K$, $L\cap N$ is reflexive.
\item[(ii)] The pair $(L,M)$ is cotauberian if and only if for every closed subspace $N$ of $Z$ with $(M,N) \in\K$, $Z/\overline{L+N}$ is reflexive.
\end{enumerate}
\end{Th}
\Proof
(i) The direct implication is a consequence of part (3) in Proposition \ref{prop:basic}: if $Q_NJ_L$ is tauberian, its kernel is reflexive.

For the converse, we assume that $(M,N)\notin \W_+$, and we will show that there exists $L$ such that $(L,M)\in\K$ and $L\cap N$ is non-reflexive.

Suppose $Q_NJ_M$ is not tauberian. By Proposition \ref{tau-op}, there exists a bounded sequence $(m_k)$ in $M$ with no weakly convergent subsequence such that $(Q_N m_k)$ converges to some $z_0+N$.
Passing to a subsequence, we can assume that $(m_k)$ is a basic sequence \cite{Alb-Kalton:06}. Hence, there exists a bounded sequence $(z^*_k)$ in $Z^*$ such that
$$
\langle m_i,z^*_j\rangle= \delta_{i,j} \textrm{ and } \|Q_N(m_k-z_0)\|\to 0.
$$
Thus, we can select a sequence $(n_k)$ in $N$ such that $\|m_k-n_k-z_0\|\to 0$. Note that $(n_k)$ is a bounded sequence in $N$ with no weakly convergent subsequence.

Passing to a subsequence, we can also assume that $\sum_{k=1}^\infty \|z^*_k\|\, \|m_k-n_k-z_0\|< 1$. Thus, $Kz= \sum_{k=1}^\infty \langle z,z^*_k\rangle (m_k-n_k-z_0)$ defines the operator $K\in \K(Z)$ with $\|K\|<1$.

By Proposition \ref{useful-prop}, $L=(I-K)(M)$ satisfies $(L,M)\in\K$. Since $(I-K)m_i=n_i-z_0$, we have $(n_i-z_0)\subset L\cap (N\oplus\langle z_0 \rangle)$. Then $L\cap (N\oplus\langle z_0\rangle)$ is non-reflexive, hence $L\cap N$ is non-reflexive.
%\textcolor{red}{$L$ appears in the statement 1}
\medskip

(ii) The direct implication is a consequence of part (4) in Proposition \ref{prop:basic}: if $Q_NJ_L$ is cotauberian, its cokernel $(Z/N)/(\overline{L+N}/N) \simeq Z/\overline{L+N}$ is reflexive.

For the converse, we assume $(L,M)\notin\W_-$ and find $N$ such that $(M,N)\in\K$ and $Z/\overline{L+N}$ is not reflexive.

We have $(M^\perp,L^\perp)\notin\W_+$, hence $Q_{L^\perp}J_{M^\perp}\notin \W_+$. So the arguments  in part (1) provide a bounded basic sequence $(m^*_k)$ in $M^\perp$ with no weakly convergent subsequence and such that $(Q_{L^\perp}m^*_k)$ converges to $z^*_0+ L^\perp\in Z^*/L^\perp$. A result of Johnson and Rosenthal \cite{JR:15} (see \cite[Lemma 3.1.19]{GM:10}) allows us to find a bounded sequence $(z_k)$ in $Z$ such that $\langle z_i, m^*_j\rangle= \delta_{i,j}$, and passing to a subsequence we can assume that $\sum_{k=1}^\infty \|z_i\|\, \|m^*_k-n^*_k-z^*_0\|<1$.

The operator $K:Z\to Z$ defined by $Kz = \sum_{k=1}^\infty \langle z, m^*_k-n^*_k-z^*_0\rangle z_k$ is compact with $\|K\|<1$, and its conjugate is given by $K^*z^*= \sum_{k=1}^\infty \langle z_k, z^*\rangle (m^*_k-n^*_k-z^*_0)$.
We have $(I-K^*)(M^\perp)=N^\perp$ for some subspace $N$ of $Z$ because $I-K^*$ is a conjugate isomorphism.
Moreover $(N^\perp,M^\perp)\in\K$ by Proposition \ref{useful-prop} and, as in the proof of (1), $L^\perp\cap N^\perp$ is non-reflexive. Hence its predual $Z/\overline{L+N}$ is not reflexive.
\fin

We  need the following fact:

\begin{lemma}\label{lemma2}
If $T\in\Lc(X,Y)$ is a bijective isomorphism and $N$ is a closed subspace of $Y$ then $(T^{-1} N)^\perp= T^*(N^\perp)$.
\end{lemma}
\Proof
Observe that $Tx \in N$ if and only if $\langle Tx,y^* \rangle =0$ for each $y^*\in N^\perp$, and this is equivalent to $x\in T^{-1}(N)$.
\fin

Next we give alternative versions of the perturbative characterizations in Theorem \ref{pert-char}.

\begin{Th}\label{pert-char-alt}
Let  $L$, $M$ and $N$ be closed subspaces of $Z$.
\begin{enumerate}
\item[(i)] The pair $(M,N)$ is not tauberian if and only if there exist a compact operator $K\in\Lc(Z)$ and a non-reflexive subspace $E\subset M$ such that $(I-K)(E) \subset N$.
\item[(ii)] The pair $(L,M)$ is not cotauberian if and only if there exist a compact operator $K\in\Lc(Z)$ and a closed subspace $F$ of $Z$ containing $M$ such that $Z/F$ is non-reflexive and $(I-K)^{-1}(F) \supset L$.
\end{enumerate}
\end{Th}
\Proof
(i) Suppose there exist a compact operator $K\in\Lc(Z)$ and a non-reflexive subspace $E\subset M$ such that $(I-K)(E) \subset N$. Take a bounded sequence $(l_k)$ in $E$ with no weakly convergent subsequence. Since $Q_N(I-K)J_E= 0$, $(Q_Nl_k)$ has a convergent subsequence. Thus $Q_NJ_M\notin \W_+$.

Conversely, assume that $(M,N)\notin \W_+$. Then there is a bounded sequence $(m_k)$ in $M$ with no weakly convergent subsequence and such that $(Q_N m_k)$ converges to $z_0+N \in Z/N$. We choose a sequence $(n_k)$ in $N$ such that $(\|m_k-z_0-n_k\|)$ converges to $0$.
%Note that the sequence $(n_k)$ is bounded and has no weakly convergent subsequence.

Passing to a subsequence, we can assume that $(m_k)$ is a basic sequence, so that there exists $(z^*_k)$ in $Z^*$ bounded such that $z^*_i(m_j)= \delta_{i,j}$, and $\sum_{k=1}^\infty \|z^*_k\|\, \|m_k-z_0-n_k\|<1$. Thus $Kz = \sum_{k=1}^\infty z^*_k(z) (m_k-z_0-n_k)$ defines a compact $K\in \Lc(Z)$ with $\|K\|<1$. In particular, $(I-K)$ is a bijective isomorphism on $Z$.

Since $(I-K)m_i = z_0+n_i$ for each $i\in \N$, $(I-K)([m_i])\subset N + \langle z_0 \rangle$. Thus $E= M\cap (I-K)^{-1}(N)$ is a non-reflexive subspace of $M$ and $(I-K)(E) \subset N$.
\medskip

(ii) Suppose there exists a compact operator $K\in\Lc(Z)$ and a closed subspace $F\supset M$ with $Z/F$ non-reflexive such that $(I-K)^{-1}(F) \supset L$. We have $F^\perp\subset M^\perp$ and, by Lemma \ref{lemma2}, $(I-K^*)(F^\perp)\subset L^\perp$. By part (1), $(M^\perp,L^\perp)$ is not tauberian, Hence $(L,M)$ is not cotauberian.

Conversely, assume that $(L,M)$ is not cotauberian. Then $(M^\perp,L^\perp)$ is not tauberian, and the arguments of the proof of part (2) of Theorem \ref{pert-char} give a bounded basic sequence $(m^*_k)$ in $M^\perp$ with no weakly convergent subsequence such that $(Q_{L^\perp}J_{M^\perp}m^*_k)$ converges to some $z^*_0+L^\perp$ in $Z^*/L^\perp$. Since $Q_{L^\perp}J_{M^\perp}$ is a conjugate operator, it takes closed balls onto closed sets. Thus $z_0^*= J_{M^\perp}m_0^*$ for some $m_0^*\in M^\perp$, and taking $(m^*_k-m^*_0)$ instead of $(m^*_k)$ we can assume $z_0^*=0$.

Moreover, passing to a subsequence, the arguments of the proof mentioned in the previous paragraph also give a bounded sequence $(z_k)$ in $Z$ such that $\langle z_i, m^*_j\rangle= \delta_{i,j}$ and $\sum_{k=1}^\infty \|z_i\|\, \|m^*_k-n^*_k\|<1$.
Thus $Kz = \sum_{k=1}^\infty \langle z, m^*_k-n^*_k \rangle z_k$ defines a compact operator $K:Z\to Z$ with $\|K\|<1$, and its conjugate is given by $K^*z^*= \sum_{k=1}^\infty \langle z_k, z^*\rangle (m^*_k-n^*_k)$.

We take $F=\bigcap_{k=1}^\infty \ker m^*_k$. Then $F \supset M$ and $Z/F$ is non-reflexive because $(m_k^*) \subset F^\perp\equiv (Z/F)^*$. Moreover $(I-K^*)m^*_j=n^*_j$ for every $j\in\N$. Since $F^\perp$ is the w$^*$-closure of the subspace generated by $(m^*_k)$, by Lemma \ref{lemma2} we get
$$
(I-K)^{-1}(F)^\perp = (I-K^*)(F^\perp) \subset L^\perp.
$$
Thus $L\subset (I-K)^{-1}(F)$, and we have showed that the subspace $F$ satisfies the required conditions.
\fin

\section{Some applications and questions}\label{sect:appl}

Here we describe some new concepts that can be defined in terms of the results of the previous sections, like notions of indecomposability of Banach spaces and closed operators in the classes $\W_+$, $\W_-$ and $\W$. We also point out some open questions.

\subsection{Incomparability of Banach spaces in terms of $\W_+$ and $\W_-$ pairs}

A Banach space $Z$ is \emph{indecomposable} if every complemented subspace of $Z$ is finite dimensional or finite codimensional.
The space $Z$ is \emph{hereditarily indecomposable} if every closed subspace of $Z$ is indecomposable, and $Z$ is \emph{quotient indecomposable} if every  quotient of $Z$ is indecomposable.
 The existence of heretitarily indecomposable and quotient indecomposable Banach spaces was proved by Gowers and Maurey \cite{GowersM:93,GowersM:97}.

It is easy to check that a Banach space $Z$ is indecomposable if and only if no pair $(M,N)\in\mathcal{F}_Z$ belongs to $\K_+\cap \K_-$. Moreover, $Z$ is hereditarily indecomposable if and only if no pair $(M,N)\in\mathcal{F}_Z$ belongs to $\K_+$, and $Z$ is quotient indecomposable if and only if no pair $(M,N)\in\mathcal{F}_Z$ belongs to $\K_-$.

We can introduce some variations of the previous notions of indecomposability by replacing the semigroups $\K_+\cap\K_+$, $\K_+$ and $\K_-$ of the classical Fredholm theory by $\W_+\cap\W_-$, $\W_+$ and $\W_-$ and other related semigroups.  Next, we give some characterizations of these variations. The proofs are not difficult so we leave them to the interested reader.

\begin{prop}
For a Banach space $Z$, the following assertions are equivalent:
\begin{enumerate}
\item[(i)] No pair of non-reflexive subspaces of $Z$ is tauberian;
\item[(ii)] If $(M,N)\in\mathcal{F}_Z$ and $M^{co}\cap N^{co}=\{0\}$ then $M^{co}=\{0\}$ or $N^{co}=\{0\}$.
\item[(iii)] If $(M,N)$ is a tauberian pair in $Z$ then $M$ or $N$ is reflexive.
%\item
\end{enumerate}
\end{prop}
%\begin{proof}  \end{proof}

Note that, since $M^{co}$ is a closed subspace of $Z^{co}$, if $M^{co}$ is dense in $Z^{co}$ then $M^{co}= Z^{co}$; or equivalently, $Z/M$ reflexive.

\begin{prop}
For a Banach space $Z$, the following assertions are equivalent:
\begin{enumerate}
\item[(i)] No pair of non-reflexive subspaces of $Z$ is cotauberian;
\item[(ii)] If $(M,N)\in\mathcal{F}_Z$ and $M^{co}+N^{co}$ is dense in $Z^{co}$ then $M^{co}=Z^{co}$ or $N^{co}=Z^{co}$.
\item[(iii)] If $(M,N)$ is a cotauberian pair in $Z$ then $Z/M$ or $Z/N$ is reflexive.
%\item
\end{enumerate}
\end{prop}
%\begin{proof} \end{proof}

Rosenthal introduced the \emph{strongly tauberian} (respectively \emph{strongly cotauberian}) operators as those $S\in\Lc(X,Y)$ for which $S^{co}$ is an isomorphism (respectively, $S^{co}$ is surjective) \cite{Rosenthal:99}.
 It is not difficult to show that a pair $(M,N)$ is strongly tauberian if and only if $M^{co} \cap N^{co}= \{0\}$ and $M^{co}+N^{co}$ is closed; and similarly, $(M,N)$ is strongly cotauberian if and only if $M^{co}+N^{co}= Z^{co}$.

\begin{prop}
Let $Z$ be a Banach space. Then no pair of non-reflexive subspaces of $Z$ is strongly tauberian if and only if given closed subspaces $M$ and $N$ of $Z$, if $M^{co}\cap N^{co}=\{0\}$ and $M^{co}+ N^{co}$ is closed then $M^{co}=\{0\}$ or $N^{co}=\{0 \}$.
\end{prop}

\begin{prop}
Let $Z$ be a Banach space. Then no pair of non-reflexive subspaces of $Z$ is strongly cotauberian if and only if given closed subspaces $M$ and $N$ of $Z$, if $M^{co}+ N^{co}=Z^{co}$ then $M^{co}=Z^{co}$ or $N^{co}=Z^{co}$.
\end{prop}

%For every weakly compactly generated Banach space $Z$ there exists a Banach space $X$ such that $X^{co}$ is isomorphic to $Z$ \cite{Bellenot:82}. This is the case when $Z$ is %reflexive or separable.

\subsection{Classes of closed operators}

Next, we introduce the classes $\W_+$, $\W_-$ and $\W$ for closed operators by means of pairs of subspaces as outlined in Section \ref{sect:prelim}.

%Given $T\in\Cl(X,Y)$, we consider the \emph{graph-norm} defined on $D(T)$ by $\|x\|_T = \|x\|+\|Tx\|$. It is easy to show that $X_T= (D(T),\|\cdot\|_T)$ is a Banach space and $T_c x=Tx$ defines an operator $T_c\in\Lc(X_T,Y)$.

\begin{defin}
Let $T\in \Cl(X,Y)$. Then $T\in \W_+$, $\W_-$ or $\W$ when the pair $(G(T),X)$ of closed subspaces of $X\times Y$ belongs to the corresponding class of pairs.
\end{defin}

The following results can be easily proved using the arguments in the previous parts of the paper:

\begin{prop}
For $T\in \Cl(X,Y)$, the following assertions are equivalent:
\begin{enumerate}
\item[(i)] $T\in \W_+$;
\item[(ii)] $G(T)^{co}\cap X^{co}= \{0\}$;
\item[(iii)] If $(x_n)\subset D(T)$ is bounded and $(Tx_n)$ is weakly convergent to $y\in Y$ then there is a subsequence of $(x_n)$ that weakly converges to $x\in D(T)$ and $Tx=y$.
%\item \textcolor{red}{Relation with $T_c$ tauberian.}
%\item
\end{enumerate}
\end{prop}
Its proof follows from \ref{tau-op} and \ref{tau-pair}. For the next propositions, use \ref{cotau-pair} 
and \ref{weakcompact-pair}.
\begin{prop}
For $T\in \Cl(X,Y)$, $T\in \W_-$ if and only if
$(Y^*,G(T)^\perp) \in \W_+$.
%\item \textcolor{red}{Relation with $T^*$ for $T\in \Cl_D(X,Y)$.}
%\item \textcolor{red}{Relation with $T_c$ cotauberian.}
%\end{enumerate}
\end{prop}

\begin{prop}
For $T\in \Cl(X,Y)$, the following assertions are equivalent:
\begin{enumerate}
\item[(i)] $T\in \W$;
\item[(ii)] $G(T)^{co}\subset X^{co}$;
\item[(iii)] if $(x_n)\subset D(T)$ and both $(x_n)$ and $(Tx_n)$ are bounded then $(x_n)$ has a weakly convergent subsequence.
%\item \textcolor{red}{Relation with $T_c$ weakly compact.}
%\item
\end{enumerate}
\end{prop}

These facts can be extended to the multivalued closed linear operators (also called closed linear relations) studied in \cite{Cross:98}.

%\subsection{$\Rr$-strictly singular and $\R$-strictly cosingular operators}
 We will denote by $\Rr$ the class of all reflexive Banach spaces
and will adopt the notation of Stephani in \cite{Stephani:94} in order to denote certain classes of operators introduced by Herman \cite{Herman:68} which are described in the following.

\begin{defin}
An operator $S\in\Lc(X,Y)$ is \emph{$\Rr$-singular,} $S\in \Rr \Sg\Sg$, if no restriction of $S$ to a non-reflexive subspace is an isomorphism.
The operator $S$ is \emph{$\Rr$-cosingular,} $S\in \Rr \Sg\Cg$, if for no non-reflexive  quotient $Y/B$ the operator  $Q_BS$ is surjective.
\end{defin}

Every weakly compact operator is $\Rr\Sg\Sg$ and $\Rr \Sg\Cg$, but it is not known whether $\Rr \Sg\Sg$ and $\Rr \Sg\Cg$ are operator ideals (stable under addition).

\subsection{Some questions}
\begin{enumerate}
%\item Are the tauberian (cotauberian) bounded operators stable under $\Rr$-$\Ss$ ($\Rr$-$\SC$) perturbations?

\item As said before, the properties of the operator classes $\Rr \Sg\Sg$ and $\Rr \Sg\Cg$ are not well known, but it is possible that the case of the respective classes of pairs be  easier. Study the $\Rr \Sg\Sg$ and the $\Rr \Sg\Cg$ pairs of subspaces.

Are the classes of pairs $\Rr \Sg\Sg$ and $\Rr \Sg\Cg$ transitive?

Is the class of tauberian pairs (respect. cotauberian pairs) stable under $\Rr \Sg\Sg$ (respect $\Rr \Sg\Cg$) perturbations?
%\end{enumerate}

%\noindent {\textbf{Further questions:}
%\begin{enumerate}
%\item We could extend the results for closed operators $T\in C(X,Y)$ to multivalued closed operators. Does this notion coincide with the one considered by Cross in \cite{Cross:92,Cross:98}?

%\item Apply the stability results for the gap: $N$ reflexive and $\delta(M,N)<1/2$ imply $M$ reflexive; $N$ super-reflexive and $\delta(M,N)<1$ imply $M$ super-reflexive.

%These facts should be useful for proving some perturbation results.

\item Find a description of the perturbation classes for the semigroups $\W_+$ and $\W_-$.

%We observe that the perturbation classes for pairs in $\K_+$ and $\K_-$ are $\Ss$ and $\SC$ respectively \cite{Gonzalez:05}, although this is not true for operators \cite{Gonzalez:03}.
\item The reference \cite{AlG:91} exhibits an operator $T\in\W_+$ such that $T^{**}\notin\W_+$, but its construction is certainly artificial and therefore, so is the corresponding pairs of subspaces derived from $T$ and $T^{**}$. Is there a natural  tauberian pair $(M,N)$ such that $(M^{\perp\perp}, N^{\perp\perp})$ is not tauberian?
\item It was proved in \cite{JohnsonNST:15} that there are tauberian operators $T:L_1(0,1)\to L_1(0,1)$ that are not upper semi-Fredholm, but the example is not explicit because the proof is probabilistic.

For $Z=L_1(0,1)$, find non-trivial examples of tauberian pairs $(M,N)\in \mathcal{F}_{L_1(0,1)}$ in $ \W_+\setminus\K_+$. A trivial example is $(M,M)$ with $M$ an infinite dimensional reflexive subspace of $L_1(0,1)$.

\item Research into  pairs of strongly tauberian and strongly cotauberian pairs would be also interesting. 
\end{enumerate}


\begin{thebibliography}{99}

%\bibitem{AG:00} P. Aiena, M. Gonz\'alez. \emph{On inessential and improjective operators.} Studia Math. 131 (1998), 271--287.

\bibitem{AGM:01} P.\ Aiena, M.\ Gonz\'alez and A.\ Mart\'\i nez-Abej\'on. \emph{Operator semigroups in Banach space theory.} Boll.\ Unione Mat.\ Ital.\ (8) 4 (2001), 157--205.

\bibitem{Alb-Kalton:06} F. Albiac, N.J. Kalton. \emph{Topics in Banach space Theory.} Springer, 2006.

%\bibitem{AAG:92} J.A. Alvarez, T. Alvarez and M. Gonz\'alez, \emph{The gap between subspaces and perturbation of non-semi-Fredholm operators,} Bull. Austral. Math. Soc. 45 (1992), 369--376.

\bibitem{AlG:91} T. Alvarez and M. Gonz\'alez. \emph{Some examples of tauberian operators.} Proc. Amer. Math. Soc. 111 (1991), 1023-1027.

%\bibitem{ArgyrosF:00} S.A. Argyros and V. Felouzis, \emph{Interpolating hereditarily indecomposable Banach spaces.} J. Amer. Math. Soc. 13(2000), 243–294.

%\bibitem{ArgyrosH:10} S.~A.\ Argyros and R.~G.\ Haydon. \emph{A hereditarily indecomposable $\mathcal{L}_\infty$-space that solves the   scalar-plus-compact problem.} Acta Math.\ 206 (2011), 1--54.

\bibitem{AsTy:90} K. Astala and H.-O. Tylli. \emph{ Seminorms related to weak compactness and to tauberian operators.} Math. Proc. Cambridge Phil. Soc. 107 (1990), 365-375.

%\bibitem{Bellenot:82} S.F. Bellenot, \emph{The $J$-sum of Banach spaces.} J.\ Funct.\ Anal.\ 48 (1982), 95--106.

%\bibitem{CarothersDil:88} N.~L.\ Carothers and S.~J.\ Dilworth. \emph{Subspaces of~$L_{p,q}$.} Proc.\ Amer.\ Math.\ Soc.\ 104 (1988), 537--545.

\bibitem{Cross:92} R.W. Cross. {\em On a theorem of Kalton and Wilansky concerning tauberian  operators.} J. Math. Anal. Appl. 171 (1992), 156--170.

\bibitem{Cross:95} R.W. Cross. {\em $F_+$-operators are tauberian.} Quaestiones Math. 18 (1995), 129--132.

\bibitem{Cross:98} R.W. Cross. {\em Multivalued linear operators.} M. Dekker, 1998.

\bibitem{DFJP} W.J. Davis, T. Figiel, W.B. Johnson, A. Pe\l czy\'nski. \emph{Factoring weakly compact operators.} J. Funct. Anal. 17 (1974), 311--327.

%\bibitem{Dilworth:90} S.~J.\ Dilworth. \emph{A scale of linear spaces related to the $L_p$ scale.} Illinois J.\ Math.\ 34 (1990), 140--158.

%\bibitem{GalegoGP:17} E.M. Galego, M. Gonz\'alez, J. Pello. \emph{On subprojectivity and superprojectivity of Banach spaces.} Results Math. 71 (2017), 1191--1205.

\bibitem{Goldberg:66} S. Goldberg. \emph{Unbounded linear operators.} McGraw-Hill, 1966.

\bibitem{Gonzalez:03} M.\ Gonz\'alez. \emph{The perturbation classes problem in Fredholm theory.} J.\ Funct.\ Anal.\ 200 (2003), 65--70.

\bibitem{Gonzalez:05} M.\ Gonz\'alez. \emph{Fredholm theory for pairs of closed subspaces of a Banach space.} J.\ Math.\ Anal.\ Appl.\ 305 (2005), 53--62.

\bibitem{GM:10} M. Gonz\'alez and A. Mart\'\i nez-Abej\'on. \emph{Tauberian operators.}
Birkh\"auser, 2010.

%\bibitem{GonzalezM-AS:10} M.\ Gonz\'alez, A.\ Mart\'\i nez-Abej\'on and M.\ Salas-Brown. \emph{Perturbation classes for semi-Fredholm operators on subprojective and superprojective spaces.} Ann.\ Acad.\ Sci.\ Fennicae Math.\ 36 (2011), 481--491.

%\bibitem{GO:89} M. Gonz\'alez and V.M. Onieva \emph{Lifting results for sequences in Banach spaces.} Math. Proc. Cambridge Phil. Soc. 105 (1989), 117-121.

\bibitem{GO:90} M. Gonz\'alez and V.M. Onieva. \emph{Characterizations of tauberian operators and other semigroups of operators.} Proc. Amer. Math. Soc. 108 (1990), 399-405.

\bibitem{GSaTy:95} M. Gonz\'alez, E. Saksman and H.-O. Tylli. \emph{Representing non-weakly compact operators.} Studia Math. 113 (1995), 265-282.

\bibitem{GowersM:93} W.~T.\ Gowers and B.\ Maurey. \emph{The unconditional basic sequence problem.} J.\ Amer.\ Math.\ Soc.\ 6 (1993), 851--874.

\bibitem{GowersM:97} W.~T.\ Gowers and B.\ Maurey. \emph{Banach spaces with small spaces of operators.} Math. Ann., 307 (1997), 543-568.

\bibitem{Herman:68} R.H. Herman. \emph{Generalizations of weakly compact operators.} Trans. Amer. Math. Soc. 132 (1968), 377--386.

\bibitem{Holub:93} J.H. Holub. {\em Characterizations of tauberian  operators and related operators on Banach spaces.} J. Math. Anal. Appl. 178 (1993), 280--288.

\bibitem{JohnsonNST:15} W.B. Johnson, A-B. Nasseri, G. Schechtman and T. Tkocz \emph{Injective Tauberian operators on $L_1$ and operators with dense range on $\ell_\infty$.} Canad. Math. Bull. 58 (2015), 276--280.

\bibitem{JR:15} W.B. Johnson and H.P. Rosenthal.  \emph{On $w^*$-basic sequences and their applications to the study of Banach spaces.} Studia Math. 43 (1972), 77--92.

%\bibitem{KO:99} N.J. Kalton and M.I. Ostrovskii. \emph{Distances between Banach spaces.} Forum Math. 11 (1999), 17--48.

\bibitem{KW:76} N. Kalton and A. Wilansky. \emph{Tauberian operators in Banach spaces.} Proc. Amer. Math. Soc. 57 (1976), 251-255.

\bibitem{Kato:66} T. Kato. \emph{Perturbation theory for linear operators.} Springer 1966.

%\bibitem{Kosz04} P. Koszmider. \emph{Banach spaces of continuous functions with few operators.} Math. Ann. 330 (2004), 151--183.

%\bibitem{LT-I} J.\ Lindenstrauss and L.\ Tzafriri. \emph{Classical Banach spaces I. Sequence spaces.} Springer, 1977.

%\bibitem{LT-II} J.\ Lindenstrauss and L.\ Tzafriri. \emph{Classical Banach spaces II. Function spaces.} Springer, 1979.

\bibitem{NR:85} R. Neidinger and H.P. Rosenthal. \emph{Norm-attainment of linear functionals on subspaces and characterizations of tauberian operators.} Pacific J. Math. 118 (1985), 215--228.

%\bibitem{O:94} M.I. Ostrovskii. \emph{Topologies on the set of all subspaces of a Banach space and related questions of Banach space geometry.} Quaestiones Math. 17 (1994), 259--319.

%\bibitem{O:96} M.I. Ostrovskii. \emph{Classes of Banach spaces stable and unstable with respect to the opening.} Quaestiones Math. 19 (1996), 191--210.

%\bibitem{O:02} M.I. Ostrovskii. \emph{Paths between Banach spaces.} Glasg. Math. J. 44 (2002), 261--273.

\bibitem{Pietsch:80} A.\ Pietsch. \emph{Operator ideals.} North-Holland, 1980.

\bibitem{Rosenthal:99} H.P. Rosenthal. \emph{On wide-(s) sequences and their applications to certain classes of operators.} Pacific J. Math. 189 (1999),  311--338.

\bibitem{Stephani:94} I. Stephani. \emph{ Operator ideals generalizing the ideal of strictly singular operators.} Math. Nachrichten 94 (1980), 29–41.


\bibitem{Tac:83} D.G. Tacon. \emph{Generalized semi-Fredholm transformations.} J. Austral. Math. Soc. A34 (1983), 60-70.

%\bibitem{Tarbard:12} M.\ Tarbard. \emph{Hereditarily indecomposable, separable $\mathcal{L}_\infty$-spaces with $\ell_1$ dual having few operators, but not very few operators.} J.\ London Math.\ Soc.\ 85 (2012), 737--764.

%\bibitem{Weis:81} L.\ Weis. \emph{Perturbation classes of semi-Fredholm operators.} Math.\ Z.\ 178 (1981), 429--442.

%\bibitem{Whitley:64} R.~J.\ Whitley. \emph{Strictly singular operators and their conjugates.} Trans.\ Amer.\ Math.\ Soc.\ 113 (1964), 252--261.

\bibitem{Yang:76} K.~W.\ Yang. \emph{The generalized Fredholm operators.} Trans.\ Amer.\ Math.\ Soc.\ 216 (1976), 313--326.

\end{thebibliography}
\end{document}